\documentclass{owrart}

\usepackage{amsmath,amsfonts}
\usepackage{graphicx}
\usepackage{colonequals}

\newtheorem{theorem}{Theorem}
\newcommand\e{\varepsilon}
\newcommand\R{\mathbb R}

\begin{document}

%% -------------------------------------------------------------------------------
%% Example abstract.
%%
%% Please use the environment talk, it has three obligatory and one
%% optional argument. The syntax is:
%% -----------------------
%% \begin{talk}[coauthors]{Name of the speaker}{Title of the talk}{Author Sorting Index}
%%      .....
%% \end{talk}
%% -----------------------
%% The names of coauthors will appear in form of :
%% (joint work with ...)
%%
%% The Author Sorting Index is only used for sorting the list of all authors
%% by their last name. It should not contain special characters like
%% accents, german umlaute, etc. (just replace \"a by a, \^a by a, etc.)
%%
%% Please use the standard thebibliography environment to include
%% your references, and try to use labels for the bibitems, which
%% are uniquely assigned to you in order to avoid conflicts with other authors.
%% You can achieve unique labels by using our on initials before every label.
%% -------------------------------------------------------------------------------

\begin{talk}[Upanshu Sharma and Manh Hong Duong~\cite{DuongPeletierSharma}]{Mark A. Peletier}
{Coarse-graining and fluctuations: Two birds with one stone}
{Mark A. Peletier}

It is well known how a duality relation of the type
\[
I^\e(\rho) = \sup_f J^\e(\rho,f)
\]
can yield the inequality $\liminf_{\e\to0} I^\e(\rho^\e)\geq I^0(\rho^0)$ when $\rho^\e$ converges to $\rho^0$ in a topology for which $J^\e(\cdot,f)$ is continuously convergent. This idea is the basis for many Gamma-convergence results. In this talk I described how this idea can be combined with the concepts of \emph{coarse-graining} and \emph{large deviations} to give a natural context in which to formulate and prove the convergence statements that constitute rigorous coarse-graining. 

We illustrate the method on a simple abstract case. 
Given a sequence of i.i.d.\ $\mathcal X$-valued stochastic Markov processes $X^{\e,i}$, indexed by $i=1,2,\dots$ and $\e>0$, we define the \emph{empirical measure} $\rho^{n,\e}$ as the $t$-parametrized curve of measures
\begin{equation}
\label{MAP:duality}
\rho^{n,\e} :[0,T]\to \mathcal P(\mathcal X), \qquad
\rho^{n,\e}_t = \frac1n \sum_{i=1}^n \delta_{X^{\e,i}_t}.
\end{equation}
For many systems of this type it has been proven that $\rho^{n,\e}$ satisfies the \emph{large-deviation principle}
\begin{equation}
\label{MAP:ldp}
\mathrm{Prob}\bigl(\rho^{n,\e}\bigr|_{[0,T]} \approx \rho^\e\bigr|_{[0,T]}\bigr) 
\sim \exp\bigl[ -nI^\e(\rho^\e)\bigr]
\qquad \text{as }n\to\infty, \text{ for fixed $\e$,}
\end{equation}
with a characterization of the rate function $I^\e$ in the form~\eqref{MAP:duality}; see e.g.~\cite{FengKurtz06}.

The rate functional $I^\e$ characterizes not only the probability of fluctuations, through~\eqref{MAP:ldp}, but also the probability-1 behaviour: this corresponds to the equation $I^\e(\rho)=0$, which  has exactly one solution, given by the equation $\partial_t \rho_t = (A^\e)^T \rho_t$. Here $A^\e$ is the generator of the processes $X^{\e,i}$.

While we describe the situation here for i.i.d.\ processes, many generalizations are available for interacting particle systems; the ideas of this talk apply to many of these systems as well.

\medskip

We define \emph{coarse-graining} as the shift to a reduced description through a \emph{coarse-graining map} $\Pi:\mathcal X\to\mathcal Y$, which typically is highly non-injective; the challenge is to characterize the behaviour of the stochastic processes $Y^{\e,i} \colonequals \Pi(X^{\e,i})$ in the limit $\e\to0$. Note that the coarse-grained equivalent of $\rho: [0,T]\to\mathcal P(\mathcal X)$ is the push-forward $\hat\rho \colonequals \Pi_\#\rho:[0,T]\to\mathcal P(\mathcal Y)$. 

The central idea of this talk is contained in the following calculation:
\begin{eqnarray*}
I^\e (\rho^\e) &=& \sup_f \;J^\e(\rho^\e,f)\\
&\stackrel{f=g\circ \Pi} \geq & \sup_g \;J^\e(\rho^\e,g\circ \Pi)\\
&\stackrel{(*)}\equalscolon& \sup_g\; {\widehat J^\e}({\hat\rho^\e},g) \\
&& \phantom{\sup\; \widehat J^\e(}\Big\downarrow\;\e\to0\\[\jot]
{\widehat I^0}({\hat \rho^0})
&\stackrel{(**)}\colonequals& 
\sup_g\; {\widehat J^0}({\hat\rho^0},g) \\
\end{eqnarray*}
The inequality above arises from the reduction to a subset of all functions $f$, namely those that are of the form $f=g\circ \Pi$. The critical step is $(*)$: here one requires that the combination of loss-of-information in passing from $\rho^\e$ to $\hat\rho^\e$ is consistent with the loss-of-resolution in considering only functions $f=g\circ \Pi$. This step essentially requires a proof of \emph{local equilibrium}; it states that the behaviour of $\rho^\e$ is such that the missing information can be deduced from the push-forward $\hat\rho^\e$, at least approximately in the limit $\e\to0$. This is at the heart of many coarse-graining methods, it is often laborious, and usually it can not be avoided.

Assuming that $\widehat J^\e(\cdot,g)$ converges in an appropriate manner to some $\widehat J^0(\cdot,g)$, we then define $\widehat I^0$ by duality in terms of $\widehat J^0$ as in $(**)$. Whether or not $\widehat I^0$ is the rate functional of some stochastic process can not be answered at this level of abstraction, and is to be determined case by case.

\bigskip

We now make the discussion more concrete by considering a specific system. Consider the stochastically perturbed Hamiltonian system
\begin{subequations}
\label{MAPpb:main}
\begin{align}
dQ &= \frac1\e P\, dt,\\
dP &= -\frac1\e \nabla V(Q)\, dt + \sqrt{2} \,dW,
\end{align}
\end{subequations}
where $P,Q$ take values in $\R$, $V\in C^2(\R)$ is a given potential with quadratic growth, and $W$ is a standard Wiener process. The stochastic differential equation~\eqref{MAPpb:main} describes a single conservative degree of freedom, such as a particle in a well or an anharmonic oscillator, with non-conservative noise; the noise appears only in the second equation, which is a force balance. 

For the purposes of illustration we will choose the double-well potential $V(q) = (q^2-1)^2/4$. With this choice the Hamiltonian $H(q,p) = p^2/2 + V(q)$ also has a double-well-structure, as is shown in Figure~\ref{MAP:fig:simulation}. 
\begin{figure}[h]
\centering
\includegraphics[height=3cm]{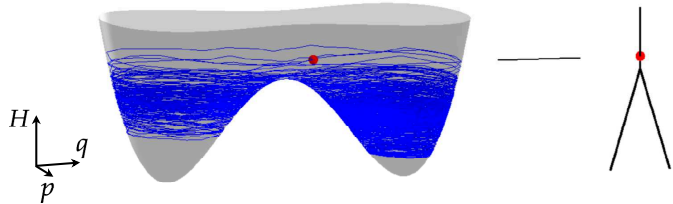}
\caption{A single run of the SDE~\eqref{MAPpb:main}. The graph on the right-hand side is a reduction of the state space $\R^2$, according to the map $\Pi$; see the text for details.}
\label{MAP:fig:simulation}
\end{figure}

Without the noise, the system is deterministic and preserves the Hamiltonian $H(q,p) = p^2/2 + V(q)$, and  solutions follow level sets of $H$. With noise, however, the Hamiltonian is not  preserved, and the solutions follow a stochastic path that stays more or less close to a level curve, depending on the size of $\e$. 

We will be interested in the limit $\e\to0$; in this limit there is a separation of time scales, in which the solutions have a fast $H$-conserving drift, and follow level sets very closely over $O(\e)$ times, with velocity $O(1/\e)$; at $O(1)$ time scales the value of~$H$ changes, and performs a biased Brownian motion, as was first proved by Freidlin and Wentzell~\cite{FreidlinWentzell94}. We  re-prove this result as an illustration of the method.

\medskip

We write $X^\e = (Q^\e,P^\e)$ for the process in $\R^2$ described by the SDE~\eqref{MAPpb:main} with a deterministic initial datum $x_0\in\R^2$, and we consider a sequence $X^{\e,i}$, $i=1,2,\dots$ of i.i.d.\ copies of this process. For this system Cattiaux and L\'eonard~\cite{CattiauxLeonard94} prove  the large-deviation principle~\eqref{MAP:ldp} and the characterization~\eqref{MAP:duality}.
\medskip

The coarse-graining map $\Pi$ maps $\R^2$ to the \emph{graph} $\Gamma$ consisting of equivalence classes of level sets of $H$, under the equivalence relation of belonging to the same connected component of the level sets of $H$. Below the saddle-point each level set has two connected components, thus leading to the two prongs in the graph $\Gamma$.
For this system the method described above yields:
\begin{theorem}
\begin{itemize}
\item $\hat\rho^\e\longrightarrow \hat \rho^0$ as $\e\to0$ in $C([0,T];\mathcal P(\Gamma))$;
\item $\liminf_{\e\to0} I^\e(\rho^\e) \geq \widehat I^0(\hat\rho^0)$.
\end{itemize}
\end{theorem}

The liminf inequality in this theorem  implies a type of convergence of solutions: the projected stochastic processes $Y^\e = \Pi(X^\e)$ converge to biased diffusions on the graph $\Gamma$, and their behvaviour is fully characterized by the law $\mu\in \mathcal P(\Gamma)$ that uniquely satisfies $\widehat I^0(\mu) = 0$. This equation can be shown to be equivalent to the diffusion-process description of~\cite{FreidlinWentzell94}, and to a weak-solution concept for a PDE. 

\bibliographystyle{plain}
\bibliography{ref}

\end{talk}

\end{document}